\documentclass[12pt]{amsart}

\usepackage{amssymb, amsthm, hyperref}
\usepackage{enumerate, amsfonts, amsxtra,amsmath,latexsym,epsfig, color}
\usepackage{mathrsfs}

\usepackage{xypic}

\input xy
\xyoption{all}

\newtheorem {theorem}{Theorem}

\newtheorem {observation} [theorem] {Observation}

\theoremstyle{definition}
\newtheorem{definition}[theorem]{Definition}
\newtheorem{remark}[theorem]{Remark}

\newtheorem{examples}[theorem]{Examples}

\def\co{\colon\thinspace}

\def\N {\mathbb N}
\def\R {\mathbb R}

\def\Mod {\mathrm{Mod}}

\begin{document}

\title[Test elements in hyperbolic groups]
{Test elements in torsion-free hyperbolic groups}

\author[Daniel Groves]{Daniel Groves}
\address{Daniel Groves\\
MSCS UIC 322 SEO, \textsc{M/C} 249\\
851 S. Morgan St.\\
Chicago, IL 60607-7045, USA}
\thanks{This work was supported by NSF grants DMS-0804365
and CAREER DMS-0953794}
\email{groves@math.uic.edu}

\begin{abstract}
We prove that in a torsion-free hyperbolic group, an element is a test element if and only if it is not contained in a proper retract.
\end{abstract}

\date{\today}

\maketitle

\begin{definition}\cite[Definition 1]{Turner}, \cite[Definition 1]{OT}
Let $G$ be a group.  An element $g \in G$ is a {\em test element} if any endomorphism $\phi \co G \to G$ for which $\phi(g) = g$ is an automorphism of $G$.
\end{definition}
This concept was studied by Shpilrain \cite{Shpilrain}, before being made explicit in \cite{Turner, OT}.  A method for constructing test elements in fee groups was given by Dold in \cite{Dold}.  Also,
Nielsen \cite{Nielsen} proved that $[a,b]$ is a test element in $F_2 = \langle a,b \mid \ \  \rangle$.  Other test elements were found by Zieschang \cite{Zieschang} and also by Shpilrain \cite{Shpilrain}.

\begin{examples}
Suppose that $F_r$ is a free group of rank $r$, with basis $\{ a_1, \ldots , a_r \}$.  For $k \ge 2$, the element $a_1^k \cdots a_r^k$ is a test element.  If $r$ is even, the element $[a_1, a_2] \cdots [a_{r-1},a_r]$ is a test element.  See \cite{Zieschang, Dold, Shpilrain, Turner}.
\end{examples}

\begin{definition}
Suppose that $G$ is a group and $H$ a subgroup, with the inclusion map $\iota \co H \to G$.  A {\em retract} is a homomorphism $r \co G \to H$ so that $r \circ \iota = {\mathrm Id}_H$.  A (proper) retract of $G$ is a (proper) subgroup $H$ for which there admits a retract $r \co G \to H$.
\end{definition}

Clearly, if $g \in G$ is contained in a proper retract of $G$, then $g$ cannot be a test element.  

\begin{definition} \cite[Definition 2]{OT}
A hyperbolic group $G$ is {\em stably hyperbolic} if for every endomorphism $\phi \co G \to G$, there are arbitrarily large values of $n$ so that $\phi^n(G)$ is hyperbolic.
\end{definition}

O'Neill and Turner \cite{OT} proved the following result.

\begin{theorem} \cite[Theorem 1]{OT} \label{t:Stable retract}
Suppose that $G$ is a torsion-free and stably hyperbolic group.  Then $g \in G$ is a test element if and only if $g$ is not contained in a proper retract of $G$.
\end{theorem}

We do not know if every torsion-free hyperbolic group is stably hyperbolic, as conjectured by O'Neill and Turner.  However, we prove that the above retract theorem holds for all torsion-free hyperbolic groups.
\begin{theorem}\label{t:SH}
Suppose that $G$ is a torsion-free hyperbolic group.  An element $g \in G$ is a test element if and only if $g$ is not contained in a proper retract of $G$.
\end{theorem}

The proof of this theorem uses Sela's Shortening Argument, and the theory of JSJ decompositions of groups.  We attempt to give references, though everything we do is standard in this area, and we assume the reader is familiar with these techniques.  For an introduction to the general theory of JSJ decompositions, see \cite{GL, GL2}.

\thanks{I would like to thank  Michael Siler, for introducing test elements to me, and for helpful discussions, and the referee for numerous useful comments and suggestions.}

\section{Proof of Theorem \ref{t:SH}}

Throughout, $G$ is a torsion-free hyperbolic group, $\phi \co G \to G$ is an endomorphism and 
$g\in G$ satisfies $\phi(g) = g$.  First note that according to the main result of \cite{Sela}, if $\phi$ is surjective then it is an automorphism.  Also, we have the following result.

\begin{theorem} \label{t:ACC} [Sela]
There exists an $N \in \N$ so that for all $n \ge N$ we have
\[	\ker(\phi^n) = \ker(\phi^N) .	\]
\end{theorem}

\begin{remark}
Theorem \ref{t:ACC} is claimed in \cite{Sela} (it does not require $\phi$ to fix any element of $G$), though a proof does not appear there.  However, if $G$ is a torsion-free hyperbolic group, then $G$ and its endomorphic images are all $G$-limit groups, in the sense of \cite[Definition 1.11]{Sela:DioHyp}.   Thus, Theorem \ref{t:ACC} is an immediate consequence of \cite[Theorem 1.12]{Sela:DioHyp}, the descending chain condition for $G$-limit groups.
\end{remark}

The sequence of kernels $\ker(\phi^i)$ is an ascending chain of subgroups of $G$.  Theorem 
\ref{t:ACC} says that this sequence stabilizes.
In particular
\[	\phi^i \left|_{\phi^N(G)} \right.	\]
is injective for all $i \ge 1$.

Consider the group $H = \phi^N(G)$, as an abstract finitely generated group.
Clearly, if we choose a different value of $N$, still satisfying the conclusion of Theorem \ref{t:ACC}, the group $H$ is unchanged (as an abstract group).

Let $\pi \co H \to \phi^N(G)$ be an isomorphism, and let $g_\pi = \pi^{-1}(g)$.

\begin{observation}
Since roots are unique in torsion-free hyperbolic groups, if $C_G(g) = \langle \gamma \rangle$, then $\phi(\gamma) = \gamma$, and $\gamma \in \phi^i(G)$ for any $i$.  Therefore, we suppose henceforth that $g$ generates its own centralizer (so that it is not a proper power in $G$).  Thus we may assume that $g_\pi$ is not a proper power in $H$.
\end{observation}

\begin{definition}
Let $\Gamma$ be a group and $\Lambda$ a subgroup of $\Gamma$.  We say that $\Gamma$ is {\em freely indecomposable rel $\Lambda$} if there is no proper free product decomposition $\Gamma = \Gamma_1 \ast \Gamma_2$ where $\Lambda \le \Gamma_1$.
\end{definition}

The relative version of Grushko's Theorem is the result below.  The proof is the same as the usual version of Grushko's Theorem, except that only free splittings where $\Lambda$ is contained in one factor are considered.  See \cite[$\S4.2$]{GL} for a discussion about why JSJ decompositions (including the Grushko decomposition) can be performed in the relative case.  The following statement can also be found in \cite{DF}.

\begin{theorem}
Let $\Gamma$ be a finitely generated group and $\Lambda$ a subgroup of $\Gamma$.  There is a free product decomposition
\[	\Gamma = \Gamma_\Lambda \ast \Gamma_1 \ast \cdots \ast \Gamma_k \ast F	\]
where 
\begin{enumerate}
\item $\Lambda \le \Gamma_\Lambda$;
\item $\Gamma_\Lambda$ is freely indecomposable rel $\Lambda$;
\item The $\Gamma_i$ are freely indecomposable and not free; and
\item $F$ is a finitely generated free group.
\end{enumerate}
The subgroup $\Gamma_\Lambda$ is unique.
Up to reordering and conjugation, the $\Gamma_i$ are unique.  The rank of $F$ is determined by $\Gamma, \Lambda$.  This splitting is called the {\em Grushko decomposition of $\Gamma$ rel $\Lambda$}.
\end{theorem}

Consider the Grushko decomposition of $H$ rel $C$, where $C = \langle g_\pi \rangle$.
The subgroup $H_C$ is freely indecomposable rel $C$ and is a retract of $H$. 

Whenever $\Gamma$ is a finitely generated group and $\Lambda$ is a subgroup, so that $\Gamma$ is freely indecomposable rel $\Lambda$, there is a relative cyclic JSJ decomposition of $\Gamma$ rel $\Lambda$.  This has the form of a graph of groups with cyclic edge groups.  There is a distinguished vertex group $V_\Lambda$, which contains $\Lambda$.  Other vertices are either cyclic, {\em QH-subgroups}, which are isomorphic to the fundamental group of a $2$-orbifold with boundary so that the adjacent edge groups correspond to boundary components or are {\em rigid} (which just means they are not of the first two types).  

That the cyclic JSJ decomposition of $H_C$ rel $C$ exists follows as in the paragraph at the end of \cite[$\S1$]{Sela:DioHyp}\footnote{This argument in turn follows that in \cite[$\S9$]{Sela:Dio1}}.  For an alternative explanation, note that since
$H_C$ is a subgroup of a torsion-free hyperbolic group, it is torsion-free and CSA.  Therefore, the existence of the required splitting follows from \cite[Theorem 11.1]{GL}.  

Let $\mathcal T(H_C,C)$ be the canonical cyclic JSJ decomposition of $H_C$ rel $C$.  Let $V_C$ be the distinguished vertex containing $C$.

The {\em modular group} of $H_C$ rel $C$, denoted $\Mod(H_C,C)$ may be defined to be the group of automorphisms of $H_C$ generated by (i) inner automorphisms of $H_C$ fixing $g_\pi$ (these are conjugation by powers of $g_\pi$); (ii) Dehn twists in edge groups of $\mathcal T(H_C,C)$; and (iii) Dehn twists in essential simple closed curves in surfaces corresponding to QH subgroups of $\mathcal T(H_C,C)$.  By convention, we choose Dehn twists which fix $V_C$ element-wise.

Suppose that $X(H_C,C) = \{ \eta \co H_C \to G \mid \eta \mbox{ injective ,} \eta(g_\pi) = g \}$.

  There is a natural action by pre-composition of $\Mod(H_C,C)$ on $X(H_C,C)$.
The Shortening Argument implies the following:

\begin{theorem} \label{t:shortening}
The set
$X(H_C,C) / \Mod(H_C,C)	$
is finite.
\end{theorem}
Theorem \ref{t:shortening} follows from the construction of the restricted Makanin-Razborov diagram for $H_C$ as in \cite[$\S1$]{Sela:DioHyp} (see also \cite[$\S8$]{Sela:Dio1} for more details in the similar situation of a free group).  This diagram encodes all of the homomorphisms from $H_C$ to $G$, where we force certain elements to have given image.  There are proper quotients of $H_C$ in this diagram, but we are only considering injective homomorphisms, so we are only concerned about the end of the diagram, which consists of finitely generated subgroups of $G$ along with injective homomorphisms into $G$.  Theorem \ref{t:shortening} is just a restatement about this last part of the restricted Makanin-Razborov diagram.

Note that normally one might expect to have to shorten by inner automorphisms of $G$, but in this case we are fixing the image of $g_\pi$, so we can only conjugate by elements centralizing $g$, and this can be achieved by inner automorphisms of $H_C$.  The limiting $\R$-tree in this construction is described in detail in the proof of Proposition 3.6 in \cite{GWilton}.

The Main Theorem is a fairly easy consequence of Theorem \ref{t:shortening}, as follows.

Suppose that $\psi_0 = \phi^N$, so that $\psi_0(G) \cong H$, and recall that $\pi \co H \to \psi_0(G)$ is an isomorphism.  Note that $\psi_0|_{\psi_0(G)}$ is injective.  Let $\eta \co H \to H_C$ be the canonical retraction and $\iota \co H_C \to H$ be the inclusion, so that $\eta \circ \iota = {\mathrm{Id}}_{H_C}$.  Let $K = \pi^{-1} (H_C)$.

We have a homomorphism $\kappa \co G \to K$ defined by
\[	\kappa = \pi^{-1} \circ \iota \circ \eta \circ \pi \circ \psi_0	.	\]
We note that $\kappa(g) = \pi^{-1} (\iota(\eta(\pi(g) )) = g$, and that 
$\kappa |_K$ is injective, since $\iota \circ \eta |_{H_C}$ is injective and $\pi$ is an isomorphism.

For a positive integer $s$, define a homomorphism $\xi_s \co H_C \to G$ by
\[	\xi_s = \kappa^s \circ \pi	.	\]
The above observations show that we have $\xi_s \in X(H_C,C)$ for any $s \ge 1$.

By Theorem \ref{t:shortening} there are positive integers $k,j$ with $k > j$ and $\alpha \in \Mod(H_C,C)$ so that 
\[	\xi_{k} = \xi_{j} \circ \alpha	.	\]

Let $\beta = \pi \circ \alpha \circ \pi^{-1}$ be the automorphism of $K$ induced by $\alpha$.
When all homomorphisms in the next equation are restricted to have $K$ as domain, we have

\[	\kappa^{k} = \xi_{k} \circ \pi^{-1} = \xi_j \circ \alpha \circ \pi^{-1} = \kappa^{j} \circ \pi \circ \alpha  \circ \pi^{-1} = \kappa^{j} \circ \beta	.	\]

Now, $\kappa$ is injective on $K$, so we have $\kappa^{k-j} |_K = \beta$, so $\kappa^{k-j}(K) = K$, and $\beta^{-1} \circ \kappa^{k-j}$ is the identity map on $K$.  

Therefore, $\beta^{-1} \circ \kappa^{k-j} \co G \to K$ is a retraction and $g \in K$.  If $\phi$ is not an automorphism then we know that it is not surjective.  Since $K \le \phi^N(G)$, in this case we clearly have $K \ne G$, so it is a proper retract.  This completes the proof of Theorem \ref{t:SH}.


\begin{thebibliography}{99}

\bibitem{DF} G. Diao and M. Feighn,  The Grushko decomposition of a finite graph of finite rank free groups: an algorithm, \textit{Geom. Topol.} {\bf 9} (2005), 1835--1880. MR2175158, Zbl 1093.20022.

\bibitem{Dold} A. Dold, Nullhomologous words in free groups which are not null homologous in any proper subgroups, \textit{Arch. Math. (Basel)} {\bf 50} (1988), 564--569. MR0948271, Zbl 0628.20026.

\bibitem{GWilton} D. Groves and H. Wilton, Conjugacy classes of solutions to equations and in equations over hyperbolic groups, \textit{J. Top.} {\bf 3} (2010), 311-332.  MR2651362, Zbl 1236.20047.

\bibitem{GL} V. Guirardel and G. Levitt, JSJ decompositions: definitions, existence, uniqueness.  I: The JSJ deformation space, preprint.  

\bibitem{GL2} V. Guirardel and G. Levitt, 	JSJ decompositions: definitions, existence, uniqueness. II: Compatibility and acylindricity, preprint.

\bibitem{Nielsen} J. Nielsen, Die Isomorphismen der allgemeinen unendlichen Gruppe mit zwei Erzeugenden, \textit{Math. Ann.} {\bf 78} (1917), 385--397.   MR1511907, JFM 46.0175.01.

\bibitem{OT} J. O'Neill and E.C. Turner, Test elements and the retract theorem in hyperbolic groups, \textit{New York J. Math.} {\bf 6} (2000), 107--117.  MR1772562, Zbl 0954.20020.

\bibitem{RS} E. Rips and Z. Sela, Structure and Rigidity in Hyperbolic Groups, I, \textit{GAFA} {\bf 4} (1994), 337--372.  MR1274119, Zbl 0818.20042.

\bibitem{Sela} Z. Sela, Endomorphisms of hyperbolic groups, I: The Hopf property, \textit{Topology} {\bf 38} (1999), 301--321. MR1660337, Zbl 0929.20033.

\bibitem{Sela:Dio1} Z. Sela, Diophantine geometry over groups I: Makanin-Razborov diagrams, \textit{Publ. Math. IHES} {\bf 93} (2001), 31--105.

\bibitem{Sela:DioHyp} Z. Sela, Diophantine geometry over groups VIII: The elementary theory of a hyperbolic group, \textit{Proc. LMS (3)} {\bf 99} (2009), 217--273.   MR2520356 ,Zbl 1241.20049.

\bibitem{Shpilrain} V. Shpilrain, Recognizing automorphisms of the free group, \textit{Arch. Math. (Basel)} {\bf 62} (1994), 385--392. MR1274742, Zbl 0802.20024.

\bibitem{Turner} E.C. Turner, Test words for automorphisms of free groups, \textit{Bull. LMS} {\bf 28} (1996), 255--263. MR1374403, Zbl 0852.20022.

\bibitem{Zieschang} H. Zieschang, \"Uber Automorphismen ebener discontinuerlicher Gruppen, \textit{Math. Ann.} {\bf 166} (1966), 148--167.  MR0201521, Zbl 0151.33102.



\end{thebibliography}
\end{document}